\documentclass{amsart}

\usepackage[utf8]{inputenc}
\usepackage[T1]{fontenc}

\usepackage{mathptmx}
\usepackage{mathtools}
\usepackage{hyperref}
\usepackage{stmaryrd}
\usepackage{amsmath}
\usepackage{amssymb}
\usepackage{amsthm}
\usepackage{tikz-cd}
\usepackage{geometry}
\usepackage{verbatim}
\newtheorem{theorem}{Theorem}[section]
\newtheorem{lemma}[theorem]{Lemma}

\newtheorem{proposition}[theorem]{Proposition}
\newtheorem{definition}[theorem]{Definition}
\newtheorem{example}[theorem]{Example}
\newtheorem{remark}[theorem]{Remark}
\newtheorem{conjecture}[theorem]{Conjecture}

\numberwithin{equation}{section}

\title{
invariant algebraic D-modules over affine algebraic groups
}

\author{Yunsong Wei}
\date{}

\begin{document}
\setcounter{section}{-1}

\begin{abstract}
We study the invariant algebraic D-modules on an affine variety under the action of an algebraic group. For linear algebraic groups with the multiplication action by themselves, such D-modules correspond to representations of their Lie algebra. For unipotent algebraic groups, we show that two invariant D-modules are isomorphic if and only if they lie in the same fiber of the GIT (Geometric Invariant Theory) quotient of the space of representations under the action of conjugation. Additionally, we classify invariant D-modules over the algebraic torus and the Borel subgroup of the general linear group. 
\end{abstract}
\maketitle
\section{Introduction}\label{SectionIntroduction}
Throughout this paper, we work over the field of complex numbers $\mathbb{C}$. Let $X$ be a smooth affine variety with a linear algebraic group $G$ acting on it. Let $D(X)$ be the algebra of differential operators on $X$. We consider $D(X)$-module structures on finitely generated free $\mathbb{C}[X]$-modules that are invariant under $G$. Details are provided in Section \ref{InvariantDmodSection}. We focus on the classification of such invariant $D(X)$ (or simply $D$)-modules up to module isomorphism. In the language of differential geometry, the module isomorphism parallels the gauge equivalence of connections.

We primarily consider $X = G$ with the action given by left multiplication. In this case, we show that an invariant $D$-module corresponds to a representation of the Lie algebra $\operatorname{Lie} G$. While two isomorphic representations induce isomorphic invariant $D$-modules, the converse is not necessarily true. If further, let $G = U$ be a unipotent group with Lie algebra $\mathfrak{n}$. We define $\operatorname{Hom}_{\text{Lie alg}}(\mathfrak{n}, \mathfrak{gl}_n(\mathbb{C}))/\!\!/\operatorname{GL}_n(\mathbb{C}) := \operatorname{Spec}(\mathbb{C}[\operatorname{Hom}_{\text{Lie alg}}(\mathfrak{n}, \mathfrak{gl}_n(\mathbb{C}))]^{\operatorname{GL}_n(\mathbb{C})})$ as the categorical quotient of $\operatorname{Hom}_{\text{Lie alg}}(\mathfrak{n}, \mathfrak{gl}_n(\mathbb{C})$) by $\operatorname{GL}_n(\mathbb{C})$, and we prove the following in Sections \ref{CategoricalQuotient} and \ref{SectionMainThm}.

\begin{theorem}
\label{mainTheorem}
Let $\mathfrak{n}$ be a complex nilpotent Lie algebra and $U$ be the corresponding unipotent group. Then 
$$\operatorname{Hom}_{\text{Lie alg}}(\mathfrak{n}, \mathfrak{gl}_n(\mathbb{C}))/\!\!/\operatorname{GL}_n(\mathbb{C}) = \{\text{Invariant D-modules of rank } n \text{ over } U\}/\sim,$$ 
where $\sim$ denotes $D(U)$-module isomorphisms. Moreover, the fiber of the quotient map 
$$\operatorname{Hom}_{\text{Lie alg}}(\mathfrak{n}, \mathfrak{gl}_n(\mathbb{C})) \twoheadrightarrow \operatorname{Hom}_{\text{Lie alg}}(\mathfrak{n}, \mathfrak{gl}_n(\mathbb{C}) )/\!\!/\operatorname{GL}_n(\mathbb{C})$$ 
is exactly an isomorphism class of invariant D-modules of rank $n$ on $U$.
\end{theorem}

Particularly, the adjoint quotient $\mathfrak{gl}_n(\mathbb{C}) \twoheadrightarrow \mathfrak{gl}_n(\mathbb{C})/\!\!/\operatorname{GL}_n(\mathbb{C})$, whose fiber consists of the matrices with the conjugate semisimple part, corresponds to an isomorphism class of invariant D-modules of rank $n$ over $\mathbb{C}$. More generally, let 
$$C_l(\mathfrak{g}) := \{(x_1, x_2, \cdots, x_l) \mid [x_i, x_j] = 0, \forall 1 \leq i < j \leq l\}$$ 
be the $l$-th commuting variety of a Lie algebra $\mathfrak{g}$. When $l=2$, this has been heavily studied by many researchers, such as \cite{gerstenhaber1961dominance}, \cite{richardson1979commuting}, \cite{popov2008irregular}, and \cite{bulois2011commuting}. Additionally, \cite{jelisiejew2022components} discusses irreducible components and reducedness for general $l$. Given \cite{jelisiejew2022components}, commuting varieties of matrices correspond to finite-dimensional modules of polynomial rings and are related to the Quot schemes. The polynomial ring is a special case of the universal enveloping algebra, while we approach this problem from the perspective of differential operators. As a corollary of Theorem \ref{mainTheorem}, we obtain 
$$C_l(\mathfrak{gl}_n(\mathbb{C}) )/\!\!/\operatorname{GL}_n(\mathbb{C}) = (\mathbb{C}^l)^n / S_n = S^n \mathbb{C}^l,$$ 
the $n$-th symmetric product of $\mathbb{C}^l$  parametrizes the isomorphism classes of invariant D-modules of rank $n$ over $\mathbb{C}^l$ with the action of translation.

The identification 
$$C_l(\mathfrak{gl}_n(\mathbb{C}) )/\!\!/\operatorname{GL}_n(\mathbb{C}) = S^n \mathbb{C}^l$$ 
is also known as the Chavelley restriction theorem for commuting schemes, as discussed in \cite{chen2021invariant}, and has been obtained in \cite{vaccarino2007linear} and \cite{10.1215/00127094-2019-0085}. More generally, we prove the following in Section \ref{CategoricalQuotient}.

\begin{proposition}
\label{quotientRepofNil}
Let $\mathfrak{n}$ be a complex nilpotent Lie algebra and let $l = \operatorname{dim} \frac{\mathfrak{n}}{[\mathfrak{n}, \mathfrak{n}]}$, then 
$$\operatorname{Hom}_{\text{Lie alg}}(\mathfrak{n}, \mathfrak{gl}_n(\mathbb{C}) )/\!\!/\operatorname{GL}_n(\mathbb{C}) = S^n \mathbb{C}^l.$$
\end{proposition}

The classification of $n$-dimensional representations of $\mathfrak{g}$ involves finding representatives for $\operatorname{Hom}_{\text{Lie alg}}(\mathfrak{g}, M_n(\mathbb{C})) / \operatorname{GL}_n(\mathbb{C})$. If $\mathfrak{g}$ is nilpotent, this classification is only complete when $\mathfrak{g} = \mathbb{C}$. In this case, we have the Jordan normal forms as representatives. For pairs of commuting matrices, finding a normal form remains an open problem. However, we consider the coarser equivalence relation induced from the isomorphism between invariant $D$-modules to simplify the issues. Besides unipotent groups, we also examine other examples in Section \ref{Moreexample}.

For \( G \) being the algebraic torus, we prove the following in Section \ref{AlgTorusSect}:

\begin{theorem}
\label{InvariantonTorus}
    Let \( T = (\mathbb{C}^*)^l \) be the algebraic torus. Then
    $$C_l(\mathrm{GL}_n(\mathbb{C}))/\mathrm{GL}_n(\mathbb{C}) = \{\text{Invariant D-modules of rank } n \text{ over } T\}/\sim,$$
    where \( \sim \) denotes the relation of \( D(T) \)-module isomorphisms. Moreover, identifying \( \operatorname{Hom}_{\text{Lie alg}}(\operatorname{Lie} T, \mathrm{GL}_n(\mathbb{C})) \) with \( C_l(\mathrm{GL}_n(\mathbb{C})/\mathrm{GL}_n(\mathbb{C}) \), the fiber of the surjective exponential map 
    $$\operatorname{Exp} : \operatorname{Hom}_{\text{Lie alg}}(\operatorname{Lie} T, \mathfrak{gl}_n(\mathbb{C})) \twoheadrightarrow \operatorname{Hom}(\operatorname{Lie} T, \mathrm{GL}_n(\mathbb{C})) / \mathrm{GL}_n(\mathbb{C})$$ 
    such that \( \operatorname{Exp}(\rho)(t) = [\operatorname{exp}(\rho(t))], t \in \operatorname{Lie} T \), is exactly an isomorphism class of invariant D-modules of rank \( n \) on \( T \).
\end{theorem}

In fact, \( C_l(\mathrm{GL}_n(\mathbb{C})/\mathrm{GL}_n(\mathbb{C}) = \operatorname{Hom}(\pi_1(T), \mathrm{GL}_n(\mathbb{C})) / \mathrm{GL}_n(\mathbb{C}) \) is the moduli space of flat connections of rank \( n \) over \( T \). Therefore, we can equivalently say that any flat connection is smoothly gauge equivalent to an invariant connection, and any two invariant connections are algebraically gauge equivalent if and only if they are smoothly gauge equivalent.

Let \( G = B \subset \operatorname{GL}_l(\mathbb{C}) \) be a Borel subgroup of the general linear group. For example, we can take \( B \) to be the group of invertible upper triangular matrices in \( \operatorname{GL}_l(\mathbb{C}) \). Let \( T = B/[B,B] = (\mathbb{C}^*)^l \) be the maximal torus of \( \operatorname{GL}_l(\mathbb{C}) \). We also compute the isomorphism classes of invariant D-modules over \( B \) in Section \ref{BorelSect} and obtain the following result:

\begin{theorem}
\label{InvOnBorel}
    Let \( B \subset \operatorname{GL}_l(\mathbb{C}) \) be a Borel subgroup and \( T = B/[B,B] \). By pulling back representations via the quotient map \( B \twoheadrightarrow T \), we have
    $$\{\text{Invariant D-modules of rank } n \text{ over } B\}/\sim = \{\text{Invariant D-modules of rank } n \text{ over } T\}/\sim.$$
\end{theorem}

In Section \ref{scsgroup}, we discuss the case when \( G \) is a simple group and prove the following:

\begin{theorem}
\label{SCSGroupDmodule}
    Let \( G \) be a connected, simply connected semisimple group. Then all invariant algebraic \( D \)-modules of rank \( n \) are isomorphic.
\end{theorem}

\begin{remark}
    The above results can be stated in terms of gauge equivalence classes of algebraic connections over the trivial \( \mathbb{G} \)-principal bundle \( G \times \mathbb{G} \rightarrow G \) when \( \mathbb{G} = \operatorname{GL}_n(\mathbb{C}) \). When \( \mathbb{G} \) is of another type of reductive group, for example, \( C_m(\operatorname{Lie}\mathbb{G}) / \!\!/_{ad} \mathbb{G} \) parametrizes \( \operatorname{Lie} \mathbb{G} \)-valued invariant flat connections over \( \mathbb{C}^m \) up to \( \mathbb{G} \)-gauge equivalence.
\end{remark}

\section{Invariant D-modules}\label{InvariantDmodSection}
Let \( X \) be a smooth affine variety. Let \( \mathbb{C}[X] \) be the algebra of regular functions on \( X \). Let \( D(X) \) be the algebra of differential operators on \( X \). Then \( \mathbb{C}[X]^{\oplus n} \) is a (left) \( D(X) \)-module naturally via differentiation, and \( D(X) \) is generated by the Lie algebra of polynomial vector fields \( \operatorname{Der}_\mathbb{C}(\mathbb{C}[X]) \) due to the smoothness of \( X \).

Let \( \Omega_{\mathbb{C}[X]} \) be the \( \mathbb{C}[X] \)-module of Kähler differentials of \( X \). Let \( M_n(\Omega_{\mathbb{C}[X]}) \) be the set of matrix-valued 1-forms on \( X \). If \( \alpha \in M_n(\Omega_{\mathbb{C}[X]}) \) satisfies \( \operatorname{d}\alpha + \alpha \wedge \alpha = 0 \), we call it a \textit{flat} 1-form. We denote the set of flat 1-forms by \( M_n(\Omega_{\mathbb{C}[X]})^f \). Given \( \alpha \in M_n(\Omega_{\mathbb{C}[X]})^f \), we can define a new \( D(X) \)-module structure on \( \mathbb{C}[X]^{\oplus n} \) extending the \( \mathbb{C}[X] \)-module structure by 
\[
v(e_1, e_2, \cdots, e_n) = (e_1, e_2, \cdots, e_n) \alpha(v),
\]
where \( e_i = (\delta_{ik})_k \in \mathbb{C}[X]^{\oplus n} \) form the standard basis of \( \mathbb{C}[X]^{\oplus n} \) and \( v \in \operatorname{Der}_\mathbb{C}(\mathbb{C}[X]) \). We denote such a new module by \( \mathbb{M}_\alpha \) and say that the rank of \( \mathbb{M}_\alpha \) is \( n \). Note that if \( \alpha = 0 \), then \( \mathbb{M}_0 = \mathbb{C}[X]^{\oplus n} \) is just the \( D(X) \)-module with the action given by differentiation.

\begin{proposition}
\label{HomCalculation}
    Let \( \alpha \in M_n(\Omega_{\mathbb{C}[X]})^f \) and \( \beta \in M_m(\Omega_{\mathbb{C}[X]})^f \). Then we have 
    $$
    \mathrm{Hom}_{D(X)}(\mathbb{M}_\alpha, \mathbb{M}_\beta) = \{\mathbb{X} \in M_{m,n}(\mathbb{C}[X]) \mid v(\mathbb{X}) = \mathbb{X} \alpha(v) - \beta(v) \mathbb{X}, \forall v \in \operatorname{Der}_\mathbb{C}(\mathbb{C}[X])\},
    $$
    where \( v(\mathbb{X}) \) is the differentiation of \( \mathbb{X} \) by \( v \) entrywise.
\end{proposition}

\begin{proof}
    Let \( \{e_1, \cdots, e_n\} \subset \mathbb{C}[X]^n \) and \( \{f_1, \cdots, f_m\} \subset \mathbb{C}[X]^m \) be the standard bases of the two free modules. Since \( \mathrm{Hom}_{D[X]}(\mathbb{M}_\alpha, \mathbb{M}_\beta) \) is a subset of \( \mathrm{Hom}_{\mathbb{C}[X]}(\mathbb{M}_\alpha, \mathbb{M}_\beta) \), we can express \( \Phi \in \mathrm{Hom}_{D[X]}(\mathbb{M}_\alpha, \mathbb{M}_\beta) \) as \( \Phi(e_1, \cdots, e_n) = (f_1, \cdots, f_m) \mathbb{X} \) for some \( \mathbb{X} \in M_{m,n} \). For \( v \in \operatorname{Der}_\mathbb{C}(\mathbb{C}[X]) \), the condition \( \Phi \circ v = v \circ \Phi \) is equivalent to $v(\mathbb{X}) = \mathbb{X} \alpha(v) - \beta(v) \mathbb{X}.    $
\end{proof}

Let \( G \) be a linear algebraic group acting on \( X \). Denote the invariant part of \( \Omega_{\mathbb{C}[X]} \) under the \( G \)-action by \( \Omega_{\mathbb{C}[X]}^G \).

\begin{definition}
    An \textit{invariant algebraic D-module} is a \( D(X) \)-module of the form \( \mathbb{M}_\alpha \) such that \( \alpha \in M_n(\Omega_{\mathbb{C}[X]}^G) \cap M_n(\Omega_{\mathbb{C}[X]})^f \), i.e., \( \alpha \) is a flat invariant matrix-valued 1-form over \( X \).
\end{definition}

Let \( \mathfrak{g} \) be the Lie algebra of \( G \). For two Lie algebras \( \mathfrak{g}_1, \mathfrak{g}_2 \), let \( \operatorname{Hom}_{\text{Lie alg}}(\mathfrak{g}_1, \mathfrak{g}_2) \) be the set of Lie algebra homomorphisms between \( \mathfrak{g}_1 \) and \( \mathfrak{g}_2 \). When \( \mathfrak{g}_2 = \mathfrak{gl}_n(\mathbb{C}) \), it consists of the \( n \)-dimensional representations of \( \mathfrak{g}_1 \). In the following, we show that invariant algebraic D-modules correspond exactly to representations of \( \mathfrak{g} \) when \( X = G \).

\begin{proposition}
\label{Rep-Form-correspondence}
    If \( X = G \) is equipped with the multiplication action, then 
    $$
    M_n(\Omega_{\mathbb{C}[G]}^G) \cap M_n(\Omega_{\mathbb{C}[G]})^f = \operatorname{Hom}_{\text{Lie alg}}(\mathfrak{g}, \mathfrak{gl}_n(\mathbb{C})).
    $$
\end{proposition}

\begin{proof}
    By definition, \( \Omega_{\mathbb{C}[G]}^G = \mathfrak{g}^* \). Let \( \{x_1, \cdots, x_l\} \) be a basis of \( \mathfrak{g} \) with structure constants \( c_{ij}^k \) given by \( [x_i, x_j] = \sum_k c_{ij}^k x_k \). Let \( \{x^1, \cdots, x^l\} \subset \mathfrak{g}^* \) be the dual basis. Since 
    \[
    \operatorname{d} x^i(x_j, x_k) = x_j(x^i(x_k)) - x_k(x^i(x_j)) - x^i([x_j, x_k]) = -c_{jk}^i,
    \]
    we have \( \operatorname{d} x^i = -\sum_{j<k} c_{jk}^i x^j \wedge x^k \).

    Now we can write \( \alpha \in M_n(\Omega_{\mathbb{C}[G]}^G) \) uniquely as \( \alpha = \sum_i A_i x^i \), where \( A_i \in M_n(\mathbb{C}) \). If \( \alpha \) is flat, then 
    \[
    \operatorname{d}\alpha + \alpha \wedge \alpha = -\sum_i A_i \sum_{j<k} c_{jk}^i x^j \wedge x^k + \sum_{j,k} A_j A_k x^j \wedge x^k = \sum_{j<k} ([A_j, A_k] - \sum_i c_{jk}^i A_i) x^j \wedge x^k = 0,
    \]
    which gives \( [A_j, A_k] = \sum_i c_{jk}^i A_i \). We can then define \( \rho \in \operatorname{Hom}_{\text{Lie alg}}(\mathfrak{g}, \mathfrak{gl}_n(\mathbb{C})) \) associated with \( \alpha \) by \( \rho(x_i) = A_i, i=1, \ldots, l \). Conversely, given \( \rho \in \operatorname{Hom}_{\text{Lie alg}}(\mathfrak{g}, M_n(\mathbb{C})) \), it is straightforward to verify that \( \alpha = \sum_i \rho(x_i) x^i \) lies in \( M_n(\Omega_{\mathbb{C}[G]}^G) \cap M_n(\Omega_{\mathbb{C}[G]})^f \).

    Since the tensor \( \sum_i x_i \otimes x^i \) does not depend on the choice of basis, the isomorphism we construct is canonical.
\end{proof}

\begin{definition}
    Let \( \rho_1, \rho_2 \in \operatorname{Hom}_{\text{Lie alg}}(\mathfrak{g}, M_n(\mathbb{C})) \) be two \( n \)-dimensional representations of \( \mathfrak{g} \). Let \( \alpha_{\rho_1}, \alpha_{\rho_2} \) be the corresponding matrix-valued 1-forms via Proposition \ref{Rep-Form-correspondence}. We say \( \rho_1, \rho_2 \) are \textit{algebraic gauge equivalent} if \( \mathbb{M}_{\alpha_{\rho_1}} \) is isomorphic to \( \mathbb{M}_{\alpha_{\rho_2}} \) as \( D(G) \)-modules.
\end{definition}

It is evident that isomorphic representations are algebraic gauge equivalent.

We also need the following theorem, which has a proof in \cite{billig2018lie}.

\begin{theorem}[\cite{billig2018lie}]
Let \( G \) be a linear algebraic group. Then the Lie algebra \( \operatorname{Der}_{\mathbb{C}}(\mathbb{C}[G]) \) of polynomial vector fields on \( G \) is a free \( \mathbb{C}[G] \)-module, generated by \( \mathfrak{g} \):
$$
\operatorname{Der}_{\mathbb{C}}(\mathbb{C}[G]) \cong \mathbb{C}[G] \otimes \mathfrak{g}.
$$
\end{theorem}

Thus, in Proposition \ref{HomCalculation}, we can relax the condition \( \forall v \in \operatorname{Der}_{\mathbb{C}}(\mathbb{C}[G]) \) to \( \forall v \in \mathfrak{g} \). Therefore, we have 
$$
\mathrm{Hom}_{D(G)}(\mathbb{M}_\alpha, \mathbb{M}_\beta) = \{\mathbb{X} \in M_{m,n}(\mathbb{C}[G]) \mid v(\mathbb{X}) = \mathbb{X} \alpha(v) - \beta(v) \mathbb{X}, \forall v \in \mathfrak{g}\}.
$$

Let \( \rho_1 \in \operatorname{Hom}_{\text{Lie alg}}(\mathfrak{g}, M_n(\mathbb{C})) \) and \( \rho_2 \in \operatorname{Hom}_{\text{Lie alg}}(\mathfrak{g}, M_m(\mathbb{C})) \) be two representations of \( \mathfrak{g} \). We also have 
$$
\mathrm{Hom}_{D(G)}(\mathbb{M}_{\alpha_{\rho_1}}, \mathbb{M}_{\alpha_{\rho_2}}) = \{\mathbb{X} \in M_{m,n}(\mathbb{C}[G]) \mid v(\mathbb{X}) = \mathbb{X} \rho_1(v) - \rho_2(v) \mathbb{X}, \forall v \in \mathfrak{g}\}.
$$

Let \( \exp: \mathfrak{g} \rightarrow G \) be the exponential map for \( G \).

\begin{definition}
    We say \( G \) is \textit{weakly exponential} if \( \exp: \mathfrak{g} \rightarrow G \) has a dense image.
\end{definition}

By Theorems A and C in \cite{djokovic1997lie}, connected simple groups and connected solvable groups are weakly exponential.

In the following proposition, we characterize the solutions for weakly exponential groups.

\begin{proposition}
\label{FindAllSol}
    Let \( G \) be a weakly exponential algebraic group. Let \( \rho_1 \in \operatorname{Hom}_{\text{Lie alg}}(\mathfrak{g}, M_n(\mathbb{C})) \) and \( \rho_2 \in \operatorname{Hom}_{\text{Lie alg}}(\mathfrak{g}, M_m(\mathbb{C})) \) be two representations of \( \mathfrak{g} \). Suppose there exist representations \( \Tilde{\rho}_1: G \rightarrow \mathrm{GL}_n(\mathbb{C}) \) and \( \Tilde{\rho}_2: G \rightarrow \mathrm{GL}_m(\mathbb{C}) \) such that \( \Tilde{\rho}_i \circ \exp = \exp \circ \rho_i \) for \( i=1,2 \). If \( \mathbb{X}: G \rightarrow M_{m,n}(\mathbb{C}) \) is an analytic matrix-valued function satisfying the following system of partial differential equations
    $$
    v(\mathbb{X}) = \mathbb{X} \rho_1(v) - \rho_2(v) \mathbb{X}, \quad \forall v \in \mathfrak{g}
    $$ 
    then
    $$
    \mathbb{X}(p) = \Tilde{\rho}_2(p^{-1}) \mathbb{X}(e) \Tilde{\rho}_1(p), \quad p \in G,
    $$
    where \( e \) is the identity of \( G \).
\end{proposition}

\begin{proof}
    It is straightforward to verify that \( \mathbb{X}(p) = \Tilde{\rho}_2(p^{-1}) \mathbb{X}(e) \Tilde{\rho}_1(p) \) satisfies the differential equations.

    For \( v \in \mathfrak{g} \), let \( \Phi(s) = \mathbb{X}(\exp(sv)), s \in \mathbb{C} \). Then 
    \[
    v(\mathbb{X})(\exp(sv)) = \frac{\operatorname{d}}{\operatorname{d}t} \bigg|_{t=0} \mathbb{X}(\exp(sv) \exp(tv)) = \frac{\operatorname{d}}{\operatorname{d}t} \bigg|_{t=0} \Phi(s+t),
    \]
    and 
    \[
    \mathbb{X}(\exp(sv)) \rho_1(v) - \rho_2(v) \mathbb{X}(\exp(sv)) = \Phi(s) \rho_1(v) - \rho_2(v) \Phi(s).
    \]

    Now the ordinary differential equation 
    \[
    \frac{\operatorname{d}}{\operatorname{d}s} \Phi(s) = \Phi(s) \rho_1(v) - \rho_2(v) \Phi(s)
    \]
    is known as a differential Sylvester equation with unique solutions \( \Phi(s) = \exp(-\rho_2(v)s) \Phi(0) \exp(\rho_1(v)s) \) (see for example \cite{behr2019solution}). Therefore, \( \mathbb{X}(\exp(sv)) = \Tilde{\rho}_2(\exp(-sv)) \mathbb{X}(e) \Tilde{\rho}_1(\exp(sv)) \). Since \( G \) is weakly exponential, we can extend \( \mathbb{X} \) to all of \( G \) by continuity.
\end{proof}

\begin{proposition}
\label{algebraic2trivial}
    Let \( \rho \in \operatorname{Hom}_{\text{Lie alg}}(\mathfrak{g}, M_n(\mathbb{C})) \) be a representation of \( \mathfrak{g} \). Suppose there exists an algebraic group homomorphism \( \Tilde{\rho}: G \rightarrow \mathrm{GL}_n(\mathbb{C}) \) such that \( \Tilde{\rho} \circ \exp = \exp \circ \rho \). Then \( \mathbb{M}_{\rho} = \mathbb{M}_{\rho_n} \), where \( \rho_n \) is the \( n \)-dimensional trivial representation of \( \mathfrak{g} \).
\end{proposition}

\begin{proof}
    By Proposition \ref{FindAllSol}, \( \mathbb{X}(p) = \Tilde{\rho}(p) \) is an algebraic isomorphism between \( \mathbb{M}_{\rho} \) and \( \mathbb{M}_{\rho_n} \).
\end{proof}

\section{Categorical quotient of the space of representations of nilpotent Lie algebra}\label{CategoricalQuotient}
In this section, we prove Proposition \ref{quotientRepofNil}.

Lie subalgebras of \( \mathfrak{gl}_n(\mathbb{C}) \) are called \textit{concrete Lie algebras}. Let \( \mathfrak{g} \subset \mathfrak{gl}_n(\mathbb{C}) \) be a concrete Lie algebra. We say \( \mathfrak{g} \) is \textit{concretely nilpotent} if every element of \( \mathfrak{g} \) is a nilpotent matrix.

\begin{lemma}
\label{nilpotencyDerived}
    Let \( \mathfrak{n} \) be a solvable Lie algebra. Let \( \rho \in \operatorname{Hom}_{\text{Lie alg}}(\mathfrak{n}, \mathfrak{gl}_n(\mathbb{C})) \) be a representation of \( \mathfrak{n} \). Then \( \rho([\mathfrak{n}, \mathfrak{n}]) \) is concretely nilpotent.
\end{lemma}

\begin{proof}
    By Lie's theorem, we can choose a basis such that every element of \( \rho(\mathfrak{n}) \) is an upper triangular matrix. Then \( \rho([\mathfrak{n}, \mathfrak{n}]) = [\rho(\mathfrak{n}), \rho(\mathfrak{n})] \) is concretely nilpotent.
\end{proof}

\begin{lemma}
\label{startLemma}
    If \( A, B \) are matrices in \( M_n(\mathbb{C}) \), let \( A = A_{\mathrm{s}} + A_{\mathrm{n}} \) be the Jordan decomposition of \( A \) with semisimple part \( A_{\mathrm{s}} \) and nilpotent part \( A_{\mathrm{n}} \). Then \( [A, [A, B]] = 0 \) implies \( [A_{\mathrm{s}}, B] = 0 \).
\end{lemma}

\begin{proof}
    Firstly, one can check that \( [A, [A, B]] = 0 \) implies \( [A^m, [A^n, B]] = 0 \) for \( m, n \in \mathbb{N} \). 

    By linearity, we have \( [p(A), [q(A), B]] = 0 \) for polynomials \( p, q \). Since \( A_{\mathrm{s}} \) is a polynomial of \( A \), we have \( [A_{\mathrm{s}}, [A_{\mathrm{s}}, B]] = 0 \). 

    Now we can assume \( A_{\mathrm{s}} \) is a diagonal matrix with entries \( A_{\mathrm{s}} = \sum_{i=1}^n \lambda_i E_{ii} \) and \( B = \sum_{j,k} b_{jk} E_{jk} \). Direct calculations show 
    \[
    [A_{\mathrm{s}}, [A_{\mathrm{s}}, B]] = \sum_{j,k} (\lambda_j - \lambda_k)^2 b_{jk} E_{jk}.
    \]

    Thus, \( (\lambda_j - \lambda_k)^2 b_{jk} = 0 \) if and only if \( (\lambda_j - \lambda_k) b_{jk} = 0 \).
\end{proof}

\begin{lemma}
\label{finiteDerivativeZeroCommute}
    If \( A, B \) are matrices in \( M_n(\mathbb{C}) \), then \( \operatorname{ad}_A^k(B) = 0 \) for some \( k > 0 \) implies \( [A_{\mathrm{s}}, B] = 0 \).
\end{lemma}

\begin{proof}
    The statement is true for \( k = 1 \) and is proved for \( k = 2 \) in Lemma \ref{startLemma}. When \( k > 2 \), \( \operatorname{ad}_A^k(B) = [A, [A, \operatorname{ad}_A^{k-2}(B)]] = 0 \) implies \( [A_{\mathrm{s}}, \operatorname{ad}_A^{k-2}(B)] = 0 = \operatorname{ad}_A^{k-2}([A_{\mathrm{s}}, B]) \), using Lemma \ref{startLemma} again. Now, by induction on \( \operatorname{ad}_A^{k-2}([A_{\mathrm{s}}, B]) = 0 \), we have \( [A_{\mathrm{s}}, [A_{\mathrm{s}}, B]] = 0 \), so finally \( [A_{\mathrm{s}}, B] = 0 \).
\end{proof} 

\begin{proposition}
\label{propOnRhos}
    Let \( \rho \in \operatorname{Hom}_{\text{Lie alg}}(\mathfrak{n}, \mathfrak{gl}_n(\mathbb{C})) \) be a representation of \( \mathfrak{n} \). We define \( \rho_{\mathrm{s}}: \mathfrak{n} \rightarrow \mathfrak{gl}_n(\mathbb{C}) \) by \( \rho_{\mathrm{s}}(x) = \rho(x)_{\mathrm{s}} \). Then
    \begin{enumerate}
        \item \label{commutinggg} \( [\rho_{\mathrm{s}}(x), \rho(y)] = 0 \) for any \( x, y \in \mathfrak{n} \).
        \item \( \rho_{\mathrm{s}} \in \operatorname{Hom}_{\text{Lie alg}}(\mathfrak{n}, \mathfrak{gl}_n(\mathbb{C})) \) and in fact \( \rho_{\mathrm{s}} \) factors through \( \frac{\mathfrak{n}}{[\mathfrak{n}, \mathfrak{n}]} \).
        \item \( \operatorname{Im} \rho_{\mathrm{s}} \) can be simultaneously diagonalizable.
    \end{enumerate}  
\end{proposition}

\begin{proof}
\begin{enumerate}
    \item Since \( \mathfrak{n} \) is nilpotent, there exists some integer \( N \) such that \( \operatorname{ad}_{\rho(x)}^N \rho(y) = 0 \). Thus, Lemma \ref{finiteDerivativeZeroCommute} implies \( [\rho_{\mathrm{s}}(x), \rho(y)] = 0 \).
    
    \item Let \( x, y \in \mathfrak{n} \). Then \( \rho_{\mathrm{s}}([x, y]) = \rho([x, y])_{\mathrm{s}} = 0 \) as \( \rho([x, y]) \) is nilpotent by Lemma \ref{nilpotencyDerived}. The result from \ref{commutinggg} implies \( [\rho_{\mathrm{s}}(x), \rho_{\mathrm{s}}(y)] = [\rho(x)_{\mathrm{s}}, \rho(y)_{\mathrm{s}}] = 0 \). Thus, \( \rho_{\mathrm{s}} \) preserves the Lie brackets. The first sentence also shows that \( \rho_{\mathrm{s}} \) factors through \( \frac{\mathfrak{n}}{[\mathfrak{n}, \mathfrak{n}]} \).

    Now we show \( \rho_{\mathrm{s}} \) is linear. Since we can assume \( \rho(x) \) and \( \rho(y) \) are upper triangular, we have \( \rho(x+y) - \rho(x)_{\mathrm{s}} - \rho(y)_{\mathrm{s}} = \rho(x)_{\mathrm{n}} + \rho(y)_{\mathrm{n}} \), which is strictly upper triangular and hence nilpotent. By \ref{commutinggg}, \( [\rho(x)_{\mathrm{s}} + \rho(y)_{\mathrm{s}}, \rho(x+y)] = 0 \). Therefore, 
    \[
    \rho_{\mathrm{s}}(x+y) = \rho(x+y)_{\mathrm{s}} = \rho(x)_{\mathrm{s}} + \rho(y)_{\mathrm{s}} = \rho_{\mathrm{s}}(x) + \rho_{\mathrm{s}}(y).
    \]

    \item The image consists of commuting diagonalizable matrices.
\end{enumerate}
\end{proof}

We borrow the following theorem from \cite{kraft2000classical}, which originated from \cite{procesi1976invariant}.

\begin{theorem}[First Fundamental Theorem for matrices (\cite{kraft2000classical})]\label{FFTmatrices}
    For every finite sequence \( i_1, i_2, \ldots, i_k \) of numbers \( 1 \leq i_\nu \leq m \), define the generalized traces 
$$
\operatorname{Tr}_{i_1 \ldots i_k}: M_n(\mathbb{C})^m \rightarrow \mathbb{C} \quad \left(A_1, \ldots, A_m\right) \mapsto \operatorname{Tr}\left(A_{i_1} A_{i_2} \cdots A_{i_k}\right).
$$
Then the ring of functions on \( M_n(\mathbb{C})^m \) which are invariant under simultaneous conjugation is generated by the invariants \( \operatorname{Tr}_{i_1 \ldots i_k} \):
$$
\mathbb{C}\left[M_n(\mathbb{C})^m\right]^{\mathrm{GL}_n(\mathbb{C})} = \mathbb{C}\left[\operatorname{Tr}_{i_1 \ldots i_k} \mid k \in \mathbb{N}, 1 \leq i_1, \ldots, i_k \leq m\right].
$$ 
\end{theorem}

Let \( \operatorname{Hom}(\mathfrak{n}, \mathfrak{gl}_n(\mathbb{C})) \) be the space of all linear maps from \( \mathfrak{n} \) to \( \mathfrak{gl}_n(\mathbb{C}) \). We have \( \operatorname{Hom}_{\text{Lie alg}}(\mathfrak{n}, \mathfrak{gl}_n(\mathbb{C})) \subset \operatorname{Hom}(\mathfrak{n}, \mathfrak{gl}_n(\mathbb{C})) \). 

Let \( \rho \in \operatorname{Hom}_{\text{Lie alg}}(\mathfrak{n}, \mathfrak{gl}_n(\mathbb{C})) \) be a representation of \( \mathfrak{n} \). Since we can assume \( \rho(x) \) is upper triangular for any \( x \in \mathfrak{n} \), we have for any \( x_1, x_2, \ldots, x_l \in \mathfrak{n} \):

$$
\operatorname{Tr}(\rho(x_1) \rho(x_2) \cdots \rho(x_l)) = \operatorname{Tr}(\rho_{\mathrm{s}}(x_1) \rho_{\mathrm{s}}(x_2) \cdots \rho_{\mathrm{s}}(x_l)).
$$

In conclusion, by Theorem \ref{FFTmatrices}, we have proved:

\begin{proposition}
\label{semisimpleSemistable}
    Let \( \rho \in \operatorname{Hom}_{\text{Lie alg}}(\mathfrak{n}, \mathfrak{gl}_n(\mathbb{C})) \) be a representation of \( \mathfrak{n} \). Then \( \rho \) and \( \rho_{\mathrm{s}} \) represent the same element in     $
    \operatorname{Hom}_{\text{Lie alg}}(\mathfrak{n}, \mathfrak{gl}_n(\mathbb{C})) / \!\!/ \operatorname{GL}_n(\mathbb{C}) \subset \operatorname{Hom}(\mathfrak{n}, \mathfrak{gl}_n(\mathbb{C})) / \!\!/ \operatorname{GL}_n(\mathbb{C}).
    $
\end{proposition}

\begin{proposition}
\label{diagonalPermute}
Let \( \rho_1, \rho_2 \in \operatorname{Hom}_{\text{Lie alg}}(\mathfrak{n}, M_n(\mathbb{C})) \) be two representations of \( \mathfrak{n} \) such that \( \rho_i = \rho_{i,\mathrm{s}} \) for \( i = 1, 2 \). Suppose for any points \( x_1, x_2, \ldots, x_k \in \mathfrak{n} \),
$$
\operatorname{Tr}(\rho_1(x_1) \rho_1(x_2) \cdots \rho_1(x_k)) = \operatorname{Tr}(\rho_2(x_1) \rho_2(x_2) \cdots \rho_2(x_k)).
$$ 
Then there exists \( g \in \mathrm{GL}_n(\mathbb{C}) \) such that \( \rho_1(x) = g \rho_2(x) g^{-1} \) for any \( x \in \mathfrak{n} \). Furthermore, if the images of \( \rho_i, i = 1, 2 \) are contained in the set of diagonal matrices, then we can choose \( g \in S_n \subset \mathrm{GL}_n(\mathbb{C}) \) to be a permutation matrix.
\end{proposition}

\begin{proof}
    We can assume the images of \( \rho_i, i = 1, 2 \) are contained in the set of diagonal matrices by Proposition \ref{propOnRhos}. Since \( \rho_i \) factors through \( \frac{\mathfrak{n}}{[\mathfrak{n}, \mathfrak{n}]} \), we can assume \( \mathfrak{n} \) is Abelian with a commuting basis \( \{x_1, \ldots, x_m\} \). Let \( (A_1, A_2, \ldots, A_m) = (\rho_1(x_1), \rho_1(x_2), \ldots, \rho_1(x_m)) \) and \( (B_1, B_2, \ldots, B_m) = (\rho_2(x_1), \rho_2(x_2), \ldots, \rho_2(x_m)) \). Then each \( B_i \) is a conjugate of \( A_i \) by a permutation matrix. We may first assume \( A_1 = B_1 \) and \( A_1 = \lambda_1 I_{n_1} \oplus \cdots \oplus \lambda_k I_{n_k} \), where \( \lambda_i \neq \lambda_j \) if \( i \neq j \). Write the corresponding \( j \)-th block of \( A_i \) (\( B_i \)) as \( A_{i}^{(j)} \) (\( B_i^{(j)} \)), \( j = 1, \ldots, k \).
    
    By assumption, $
    \operatorname{Tr}(A_1^{l_1} A_2^{l_2} \cdots A_m^{l_m}) = \operatorname{Tr}(A_1^{l_1} B_2^{l_2} \cdots B_m^{l_m})
    $    holds for each \( l_1 \). 
    Then $
    \operatorname{Tr}(A_1^{l_1} A_2^{l_2} \cdots A_m^{l_m}) = \sum_{j=1}^k \lambda_j^{l_1} \operatorname{Tr}((A_2^{(j)})^{l_2} \cdots (A_m^{(j)})^{l_m}) = \sum_{j=1}^k \lambda_j^{l_1} \operatorname{Tr}((B_2^{(j)})^{l_2} \cdots (B_m^{(j)})^{l_m}) = \operatorname{Tr}(A_1^{l_1} B_2^{l_2} \cdots B_m^{l_m}), $
    for all \( j \). Thus, \( \operatorname{Tr}((A_2^{(j)})^{l_2} \cdots (A_m^{(j)})^{l_m}) = \operatorname{Tr}((B_2^{(j)})^{l_2} \cdots (B_m^{(j)})^{l_m}) \) holds for all \( l_2, \ldots, l_m \).

    Using the fundamental theorem for matrices and induction, this implies the \( j \)-th block tuples \( (A_2^{(j)}, \ldots, A_m^{(j)}) = g_j(B_2^{(j)}, \ldots, B_m^{(j)}) g_j^{-1} \) for some \( g_j \in S_{n_j} \subset GL_{n_j}(\mathbb{C}) \). Then \( g = g_1 \oplus \cdots \oplus g_k \) is the permutation required.
\end{proof}
\begin{proof}[Proof of Proposition \ref{quotientRepofNil}]
    Let $\rho_1, \rho_2 \in \operatorname{Hom}_{\text{Lie alg}}(\mathfrak{n}, M_n(\mathbb{C}))$ be two representations of $\mathfrak{n}$. Combining Proposition \ref{propOnRhos}, Proposition \ref{semisimpleSemistable}, and Proposition \ref{diagonalPermute}, we have that $\rho_1$ and $\rho_2$ represent the same element in $
    \operatorname{Hom}_{\text{Lie alg}}(\mathfrak{n}, \mathfrak{gl}_n(\mathbb{C})) / \!\!/ \operatorname{GL}_n(\mathbb{C})
    $
    if and only if there exists $g \in \mathrm{GL}_n(\mathbb{C})$ such that $\rho_{1,s}(x) = g \rho_{2,s}(x) g^{-1}$ for any $x \in \mathfrak{n}$. 

    Now notice that 
    $
    \operatorname{Hom}_{\text{Lie alg}}\left(\frac{\mathfrak{n}}{[\mathfrak{n}, \mathfrak{n}]}, \mathfrak{gl}_n(\mathbb{C})\right) = C_l(\mathfrak{gl}_n(\mathbb{C}))
    $
    after choosing a basis of $\mathfrak{n}$ and $(\mathbb{C}^l)^n = (\mathbb{C}^n)^l$ such that $\mathbb{C}^n$ represents the set of diagonal matrices in $\mathfrak{gl}_n(\mathbb{C})$.
\end{proof}

\section{Gauge equivalence classes of representations of nilpotent Lie algebra.}\label{SectionMainThm}
In this section, we prove Theorem \ref{mainTheorem}.

Let $U$ be a unipotent group with corresponding nilpotent Lie algebra $\mathfrak{n}$. Notice that the exponential map $\operatorname{exp}:\mathfrak{n}\rightarrow U$ is an algebraic isomorphism. By pull-back of functions through $\operatorname{exp}^*:\mathbb{C}[U]\rightarrow \mathbb{C}[\mathfrak{n}]$, we can reformulate the results in Section \ref{InvariantDmodSection} in terms of $\mathfrak{n}$ only.

Let $\phi\in\mathbb{C}[\mathfrak{n}]$ and $v\in\mathfrak{n}$. The differentiation $v(\phi)\in\mathbb{C}[\mathfrak{n}]$ of $\phi$ by $v$ is given by 
\begin{equation}
\label{Diff}
  v(\phi)(p) := \frac{\operatorname{d}}{\operatorname{d}t}\bigg|_{t=0} \phi(\log(\exp(p)\exp(tv))), \quad p\in \mathfrak{n}. 
\end{equation}

Note that $\log(\exp{p}\exp{tx})$ can also be defined directly by the Baker–Campbell–Hausdorff formula, as there are only finitely many terms involving Lie brackets due to the nilpotency of $\mathfrak{n}$. 

Now, extending the differentiation to matrix-valued functions, we can reformulate Proposition \ref{HomCalculation} for the unipotent group $U$. 

\begin{proposition}
\label{HomCalculationOfNil}
     Let $\alpha\in M_n(\Omega_{\mathbb{C}[U]})^f$ and $\beta\in M_{m}(\Omega_{\mathbb{C}[U]})^f$. Then we have 
    $$
    \mathrm{Hom}_{D(U)}(\mathbb{M}_\alpha,\mathbb{M}_\beta) = \{\mathbb{X}\in M_{m,n}(\mathbb{C}[\mathfrak{n}]) \mid v(\mathbb{X}) = \mathbb{X}\alpha(v) - \beta(v)\mathbb{X}, \forall v\in \mathfrak{n}\},
    $$
    where $v(\mathbb{X})$ denotes $\mathbb{X}$ differentiated by $v$ entrywise through the formula \eqref{Diff}.
\end{proposition}

We also reformulate Proposition \ref{FindAllSol} in this particular case.

\begin{lemma}
    Let $\rho_1, \rho_2 \in \operatorname{Hom}_{\text{Lie alg}}(\mathfrak{n}, M_n(\mathbb{C}))$ be two representations of $\mathfrak{n}$. If $\mathbb{X}:\mathfrak{n}\rightarrow M_n(\mathbb{C})$ is an analytic matrix-valued function satisfying the following system of partial differential equations
    $$
    v(\mathbb{X}) = \mathbb{X}\rho_1(v) - \rho_2(v)\mathbb{X}, \quad \forall v \in \mathfrak{n},
    $$
    then 
    $$
    \mathbb{X}(p) = \exp(-\rho_2(p))\mathbb{X}(0)\exp(\rho_1(p)), \quad p\in\mathfrak{n}.
    $$
\end{lemma}

\begin{proof}
     One just needs to notice that $\rho(\log(\exp(x)\exp(y)))=\log(\exp(\rho(x))\exp(\rho(y)))$ for $x,y\in\mathfrak{n}$ and $\rho\in\operatorname{Hom}_{\text{Lie alg}}(\mathfrak{n},M_n(\mathbb{C}))$, and prove this word for word by copying from Proposition \ref{FindAllSol}.     
\end{proof}

Now we can classify invariant D-modules over $U$.

\begin{proposition}
\label{semisimpleGauge}
   Let $\rho\in \operatorname{Hom}_{\text{Lie alg}}(\mathfrak{n}, \mathfrak{gl}_n(\mathbb{C}))$ be a representation of $\mathfrak{n}$. Then $\rho$ and $\rho_{\mathrm{s}}$ are algebraic gauge equivalent.
\end{proposition}

\begin{proof}
    Let $\mathbb{X}(p) = \exp(-\rho_{\mathrm{s}}(p))\exp(\rho(p)) = \exp(\rho(p)_{\mathrm{n}}), p\in\mathfrak{n}$. We claim that the map $p \mapsto \exp(\rho(p)_{\mathrm{n}})$ is algebraic, so it could serve as the isomorphism between $\mathbb{M}_{\alpha_\rho}$ and $\mathbb{M}_{\alpha_{\rho_{\mathrm{s}}}}$.

    Let $\{x_1, \ldots, x_l\}$ be a basis of $\mathfrak{n}$. Suppose $p = \sum_i a_i x_i$. We have $[\rho(p), \rho(p) - \sum_i a_i \rho(x_i)_{\mathrm{s}}] = 0$ by Proposition \ref{propOnRhos}. We can assume the image of $\rho$ lies in the space of upper triangular matrices by Lie's theorem, so $\rho(p) - \sum_i a_i \rho(x_i)_{\mathrm{s}} = \sum_i a_i \rho(x_i)_{\mathrm{n}}$ is nilpotent. Thus, $\rho(p)_{\mathrm{n}} = \sum_i a_i \rho(x_i)_{\mathrm{n}}$. 

    Now we have 
    $$
    \exp(\rho(p)_{\mathrm{n}}) = \sum_{k=0}^n \frac{(\rho(p)_{\mathrm{n}})^k}{k!} = \sum_{k=0}^n \frac{(\sum_{i=1}^l a_i \rho(x_i)_{\mathrm{n}})^k}{k!}
    $$ 
    which is polynomial in $a_1, a_2, \ldots, a_l$ and hence algebraic. 
\end{proof}

\begin{proposition}
\label{twodiagonalGauge}
    Let $\rho_1, \rho_2 \in \operatorname{Hom}_{\text{Lie alg}}(\mathfrak{n}, M_n(\mathbb{C}))$ be two representations of $\mathfrak{n}$ such that $\rho_i = \rho_{i,\mathrm{s}}$ for $i=1,2$. Suppose $\rho_1$ and $\rho_2$ are gauge equivalent. Then there exists $g \in \mathrm{GL}_n(\mathbb{C})$ such that $\rho_1(x) = g \rho_2(x) g^{-1}$ for any $x \in \mathfrak{n}$. Furthermore, if the images of $\rho_i, i=1,2$ are contained in the set of diagonal matrices, then we can choose $g \in S_n \subset \mathrm{GL}_n(\mathbb{C})$ to be a permutation matrix.
\end{proposition}

\begin{proof}
    We can assume that the images of $\rho_i, i=1,2$ are contained in the set of diagonal matrices by conjugation and Proposition \ref{propOnRhos}. We can write $\rho_i(p) = \operatorname{diag}(\lambda_1^{(i)}(p), \lambda_2^{(i)}(p), \ldots, \lambda_n^{(i)}(p))$, $i=1,2$, for $p \in \mathfrak{n}$ as diagonal matrices with one-dimensional representations $\lambda_k^{(i)}$ of $\mathfrak{n}$.

    An isomorphism between $\mathbb{M}_{\alpha_{\rho_1}}$ and $\mathbb{M}_{\alpha_{\rho_2}}$ is given by some $\mathbb{X}: p \mapsto \exp(-\rho_2(p))\mathbb{X}(0)\exp(\rho_1(p)), p \in \mathfrak{n}$ for some $\mathbb{X}(0) = (u_{ij})_{ij} \in \operatorname{GL}_n(\mathbb{C})$. Since $\operatorname{det}(\mathbb{X}(0)) \neq 0$, we have for some permutation $\sigma \in S_n$, $\prod_{i=1}^n u_{i\sigma(i)} \neq 0$, hence $u_{i\sigma(i)} \neq 0$. Now, 
    $$
    \exp(-\rho_2(p))\mathbb{X}(0)\exp(\rho_1(p)) = (u_{ij} \exp(\lambda_j^{(1)}(p) - \lambda_i^{(2)}(p))) \in M_n(\mathbb{C}[\mathfrak{n}]),
    $$ 
    we must have $\lambda_{\sigma(i)}^{(1)}(p) - \lambda_i^{(2)}(p) = 0$. Therefore, we can choose $g = \sigma^{-1}$ to be the corresponding permutation matrix.
\end{proof}

\begin{proof}[Proof of Theorem \ref{mainTheorem}]
    Comparing Proposition \ref{semisimpleSemistable} with Proposition \ref{semisimpleGauge} and Proposition \ref{diagonalPermute} with Proposition \ref{twodiagonalGauge}, we have completed the proof.
\end{proof}

\begin{remark}
  Note that by pull-back of 1-forms via the quotient $U \rightarrow U/[U,U]$ and pull-back of representations via $\mathfrak{n} \rightarrow \mathfrak{n}/[\mathfrak{n},\mathfrak{n}]$, the identification in Theorem \ref{mainTheorem} is compatible in the sense of the following diagram:
$$
\begin{tikzcd}
    \operatorname{Hom}_{\text{Lie alg}}\left(\frac{\mathfrak{n}}{[\mathfrak{n},\mathfrak{n}]}, \mathfrak{gl}_n(\mathbb{C})\right)/\!\!/ \operatorname{GL}_n(\mathbb{C}) 
    \arrow[d, "\simeq"] 
    \arrow[r, "\cong"] 
    & 
    \{\text{Invariant D-modules of rank } n \text{ over } U/[U,U] \}/\sim 
    \arrow[d, "\simeq"] \\ 
    \operatorname{Hom}_{\text{Lie alg}}(\mathfrak{n}, \mathfrak{gl}_n(\mathbb{C})) / \!\!/ \operatorname{GL}_n(\mathbb{C}) 
    \arrow[r, "\cong"] 
    & 
    \{\text{Invariant D-modules of rank } n \text{ over } U \}/\sim 
\end{tikzcd}
$$
  
\end{remark}
\section{More examples}\label{Moreexample}
\subsection{On algebraic torus}\label{AlgTorusSect}
In this subsection, we prove Theorem \ref{InvariantonTorus}.

Let $T = (\mathbb{C}^*)^l$ be an algebraic torus of dimension $l$. Let $\rho_1, \rho_2 \in \operatorname{Hom}_{\text{Lie alg}}(\operatorname{Lie}(T), M_n(\mathbb{C}))$ be two representations of $\operatorname{Lie}(T)$.

If $\rho_1$ and $\rho_2$ are algebraic gauge equivalent, then the two connections $\nabla_1 = d + \alpha_{\rho_1}$ and $\nabla_2 = d + \alpha_{\rho_2}$ must be smoothly gauge equivalent. By the fact that 
$
C_l(\mathrm{GL}_n(\mathbb{C})) / \mathrm{GL}_n(\mathbb{C}) = \operatorname{Hom}(\pi_1(T), \mathrm{GL}_n(\mathbb{C})) / \mathrm{GL}_n(\mathbb{C})
$
is the moduli space of flat connections of rank $n$ over $T$, there exists an element $g \in \mathrm{GL}_n(\mathbb{C})$ such that 
$
\exp(2\pi i \rho_1(t)) = g \exp(2\pi i \rho_2(t)) g^{-1}, \quad \forall t \in \operatorname{Lie} T.
$ 
The converse is also true.

\begin{proposition}
    Let $\rho_1, \rho_2 \in \operatorname{Hom}_{\text{Lie alg}}(\operatorname{Lie}(T), M_n(\mathbb{C}))$ be two representations of $\operatorname{Lie}(T)$. If there exists an element $g \in \mathrm{GL}_n(\mathbb{C})$ such that 
    $$
    \exp(2\pi i \rho_1(t)) = g \exp(2\pi i \rho_2(t)) g^{-1}, \quad \forall t \in \operatorname{Lie} T,
    $$ 
    then $\rho_1$ and $\rho_2$ are algebraic gauge equivalent.
\end{proposition}

\begin{proof}
    We may assume $g = I_n$ by conjugating $\rho_2$. Let $e_1, \ldots, e_l$ be a basis of $\operatorname{Lie}(T)$ and let $\rho_1(e_i) = A_i, \rho_2(e_i) = B_i$ for $i = 1, 2, \ldots, l$. Under this basis, we write $T = (\mathbb{C}^*)^l$ and let $z = (z_1, \ldots, z_n) \in T$. Let $\operatorname{Log}:T\rightarrow \operatorname{Lie}(T)$ be the multi-valued logarithm map.    

    Let $\log(z) = \sum_i \log(z_i)e_i$ be a chosen particular value in $\operatorname{Log}(z)$ for each $z \in T$. We define the map 
    $$\mathbb{X}: z \mapsto \exp(-\rho_2(\log(z)))\exp(\rho_1(\log(z))) = \exp\left(-\sum_i \log(z_i) B_i\right) \exp\left(\sum_i \log(z_i) A_i\right), \quad z \in T.$$
    In the following, we will show that $\mathbb{X} \in \mathrm{GL}_n(\mathbb{C}[T])$.

    According to the condition, we have $\exp(2\pi i A_i) = \exp(2\pi i B_i)$. By the Jordan-Chevalley decomposition, we have $\exp(2\pi i A_{i,s}) = \exp(2\pi i B_{i,s})$ and $\exp(2\pi i A_{i,n}) = \exp(2\pi i B_{i,n})$. By the Springer isomorphism between nilpotent varieties and unipotent varieties (see \cite{springer1969unipotent}, or see \cite{mcninch2009nilpotent} and \cite{sobaje2015springer}), we have $A_{i,n} = B_{i,n}$. Therefore, we can assume $A_i = A_{i,s}$ and $B_i = B_{i,s}$ for simplicity.
    
    Since $A_1, \ldots, A_l$ and $B_1, \ldots, B_l$ are two tuples of commuting diagonalizable matrices, we can assume $A_1, \ldots, A_l$ are diagonal and there exists $h \in \mathrm{GL}_n(\mathbb{C})$ such that $h B_1 h^{-1}, \ldots, h B_l h^{-1}$ are also diagonal. Write $A_i = \operatorname{diag}(\lambda^{(i)}_1, \ldots, \lambda^{(i)}_n)$ and $h B_i h^{-1} = \operatorname{diag}(\mu^{(i)}_1, \ldots, \mu^{(i)}_n)$ for $i=1, \ldots, l$. By $\exp(2\pi i A_i) = \exp(2\pi i B_i)$, we have $h_{jk} = \exp(2\pi i (\mu_j^{(i)} - \lambda_k^{(i)})) h_{jk}$, which implies that $\mu_j^{(i)} - \lambda_k^{(i)} \in \mathbb{Z}$ if $h_{jk} \neq 0$. 
    
    Now we can express 
    $$
    \begin{aligned}
    \exp(-\log(z_i) B_i) \exp(\log(z_i) A_i) &= h^{-1} \exp(-\log(z_i) \operatorname{diag}(\mu^{(i)}_1, \ldots, \mu^{(i)}_n)) h \exp(\log(z_i) \operatorname{diag}(\lambda^{(i)}_1, \ldots, \lambda^{(i)}_n)) \\
    &= h^{-1} (z_i^{\mu_j^{(i)} - \lambda_k^{(i)}} h_{jk})_{jk}  
    \end{aligned}
    $$
    which is polynomial in $z_i^\pm$. We can arrange by commutativity that 
    $$
    \mathbb{X}(z) = \exp\left(-\sum_{j \neq i} \log(z_j) B_j\right) \exp\left(-\log(z_i) B_i\right) \exp\left(\log(z_i) A_i\right) \exp\left(\sum_{j \neq i} \log(z_j) A_j\right)
    $$ 
    shows that $\mathbb{X}(z) = \exp\left(-\sum_{j \neq i} \log(z_j) B_j\right) h^{-1} (z_i^{\mu_j^{(i)} - \lambda_k^{(i)}} h_{jk})_{jk} \exp\left(\sum_{j \neq i} \log(z_j) A_j\right)$ is polynomial in each $z_i$ and $z_i^{-1}$, $i=1,\ldots,l$, with bounded degrees independent of $z_j$ for $j \neq i$. Thus, $\mathbb{X} \in \mathrm{GL}_n(\mathbb{C}[T])$.  
    
    Now it is straightforward to verify that 
    $$
    v(\mathbb{X}) = \mathbb{X} \rho_1(v) - \rho_2(v) \mathbb{X}, \quad \forall v \in \operatorname{Lie}(T)
    $$ 
    since $\mathbb{X}$ is algebraic.
\end{proof}

To finish the proof of Theorem \ref{InvariantonTorus}, it remains to clarify that the exponential map 
$$
\operatorname{Exp}: \operatorname{Hom}_{\text{Lie alg}}(\operatorname{Lie}(T), \mathfrak{gl}_n(\mathbb{C})) \twoheadrightarrow \operatorname{Hom}(\operatorname{Lie}(T), \mathrm{GL}_n(\mathbb{C})) / \mathrm{GL}_n(\mathbb{C})
$$ 
is surjective.

\begin{proposition}
    The multi-exponential map 
    $$
    \begin{aligned}
        \operatorname{Exp}: C_l(M_n(\mathbb{C})) & \rightarrow C_l(\mathrm{GL}_n(\mathbb{C})) \\
        (A_1, \ldots, A_l) & \mapsto (\exp(A_1), \ldots, \exp(A_l))
    \end{aligned}
    $$ 
    is surjective.
\end{proposition}

\begin{proof}
    We use induction on $n$ to find a preimage of $(B_1, \ldots, B_l) \in C_l(\mathrm{GL}_n(\mathbb{C})$. If $n=1$, then this is certainly correct. If $B_1$ has distinct eigenvalues, since commuting matrices preserve generalized eigenspaces, we can decompose the other $B_i$'s and then use induction. Now we may assume $B_i = \lambda_i I_n + N_i$, where $\lambda_i \neq 0$, $N_i$ is a nilpotent matrix, and $[N_i, N_j] = 0$ for $1 \leq i,j \leq l$.    
  
    Let $\log_n(x) := \sum_{i=1}^n (-1)^{i-1} \frac{(x-1)^i}{i}$, then $[\log_n(I_n + N_i/\lambda_i), \log_n(I_n + N_j/\lambda_j)] = 0$ for all $1 \leq i,j \leq l$. Therefore, 
    $$
    \left(\log(\lambda_1)I_n + \log_n(I_n + N_1/\lambda_1), \ldots, \log(\lambda_l)I_n + \log_n(I_n + N_l/\lambda_l)\right) \in \operatorname{Exp}^{-1}(B_1, \ldots, B_l).
    $$
\end{proof}
\subsection{On Borel subgroup of general linear group}\label{BorelSect}
In this subsection, we prove Theorem \ref{InvOnBorel}.

Let $E_{ij} = (\delta_{is} \delta_{jr})_{rs} \in M_l(\mathbb{C})$. Write 
$$
B = \left\{ \sum_{1 \leq i \leq j \leq l} x_{ij} E_{ij} \mid x_{ii} \neq 0, x_{ij} \in \mathbb{C}, 1 \leq i \leq j \leq l \right\}
$$ 

and the Lie algebra 
$$
\operatorname{Lie}(B) = \langle E_{kl}, 1 \leq k \leq l \leq l \rangle_{\mathbb{C}}.
$$ 

Let 
$$
T = \left\{ \sum_{i=1}^l x_{ii} E_{ii} \mid x_{ii} \in \mathbb{C}^*, 1 \leq i \leq l \right\} = B / [B, B].
$$

Let $\rho \in \operatorname{Hom}_{\text{Lie alg}}(\operatorname{Lie}(B), M_n(\mathbb{C}))$ be a representation of $\operatorname{Lie}(B)$. Define a new representation $\rho_a$ such that $\rho_a(E_{ij}) = 0$ if $i < j$ and $\rho_a(E_{ii}) = \rho(E_{ii})$. We want to show that $\rho$ is algebraic gauge equivalent to $\rho_a$. For simplicity, we write $\rho(E_{ij})$ as $A_{ij}$.

Let $\mathbb{X}: B \rightarrow \mathrm{GL}_n(\mathbb{C})$ be defined as 
$$
    \begin{aligned}
        &\mathbb{X}(x_{11},x_{12},x_{22},\cdots,\cdots,x_{1k},x_{2k},\cdots,x_{kk},\cdots,\cdots,x_{1l},x_{2l},\cdots,x_{ll})
    \\=&(\prod_{i=1}^l\exp(-{\log(x_{ii})}A_{ii}))[\exp({\log(x_{ll})}A_{ll})(\prod_{i=1}^{l-1}\exp(A_{il}x_{il}))\cdots\exp({\log(x_{kk})}A_{kk})(\prod_{i=1}^{k-1}\exp(A_{ik}x_{ik}))\cdots\\ &\cdots\exp({\log(x_{22})}A_{22})\exp(A_{12}x_{12}))\exp({\log(x_{11})}A_{11})]
    \end{aligned}
    $$

\begin{proposition}
    \label{algebraicX}
    We have $\mathbb{X} \in \mathrm{GL}_n(\mathbb{C}[B])$.
\end{proposition}

We need the following well-known lemma:

\begin{lemma}
\label{ABBlemma}
    Let $X, Y \in M_n(\mathbb{C})$, and $[X, Y] = sY$, where $s$ is a complex number. Then 
    $$
    \exp(X) \exp(Y) \exp(-X) = \exp(\exp(s) Y).
    $$
\end{lemma}

\begin{proof}
    Since $XY = Y(X + sI)$, we have $X^n Y^m = Y^m (X + msI)^n$. Now take the Taylor expansion of the exponential functions.
\end{proof}

\begin{proof}[Proof of Proposition \ref{algebraicX}]
    We can conjugate every term, so it suffices to show that $\exp(-\log(x_{ii}) A_{ii}) \exp(\log(x_{ik}) A_{ik}) \exp(\log(x_{ii}) A_{ii})$ is algebraic in $x_{ii}$ and $x_{ik}$ for $1 \leq i < k \leq l$. By Lemma \ref{ABBlemma}, we have 
    $$
    \exp(-\log(x_{ii}) A_{ii}) \exp(x_{ik} A_{ik}) \exp(\log(x_{ii}) A_{ii}) = \exp\left(\frac{x_{ik}}{x_{ii}} A_{ik}\right).
    $$
    Since $B$ is solvable, by Lemma \ref{nilpotencyDerived}, $A_{ik}$ is nilpotent, and hence $\exp\left(\frac{x_{ik}}{x_{ii}} A_{ik}\right)$ is polynomial in $\frac{x_{ik}}{x_{ii}}$.
\end{proof}

Finally, it is straightforward to verify that 
$$
v(\mathbb{X}) = \mathbb{X} \rho(v) - \rho_a(v) \mathbb{X}, \quad v \in \operatorname{Lie}(B).
$$ 
Hence, $\rho$ is algebraic gauge equivalent to $\rho_a$ by the isomorphism $\mathbb{X}$. Now, $\rho_a$ can be viewed as the representation induced from $T$ through pull-back or extension. Since the restriction of algebraic gauge equivalent representations will still be algebraic gauge equivalent, we finish the proof of Theorem \ref{InvOnBorel}.
\subsection{On simply-connected semisimple groups}\label{scsgroup}

Let $G$ be a connected semisimple group. For example, $G = \mathrm{SL}_l(\mathbb{C}), \mathrm{Sp}_{2l}(\mathbb{C}), \mathrm{SO}_l(\mathbb{C})$.

\begin{proof}[Proof of Theorem \ref{SCSGroupDmodule}]
    If $G$ is simply connected, then any representation $\rho \in \operatorname{Hom}_{\text{Lie alg}}(\mathfrak{g}, M_n(\mathbb{C}))$ integrates to give an algebraic group homomorphism $\Tilde{\rho}: G \rightarrow \mathrm{GL}_n(\mathbb{C})$ such that $\Tilde{\rho} \circ \exp = \exp \circ \rho$. Now use Corollary \ref{algebraic2trivial}.
\end{proof}

\begin{example}
    Let $\mathfrak{g} = \mathfrak{sl}_2(\mathbb{C})$ with basis 
    $$
    X = \begin{bmatrix}
        0 & 1 \\ 
        0 & 0 
    \end{bmatrix}, \quad 
    Y = \begin{bmatrix}
        0 & 0 \\ 
        1 & 0 
    \end{bmatrix}, \quad 
    H = \begin{bmatrix}
        1 & 0 \\ 
        0 & -1 
    \end{bmatrix}.
    $$ 

    Suppose $G = \mathrm{SL}_2(\mathbb{C}) = \left\{ \begin{bmatrix}
        x & y \\ 
        z & w 
    \end{bmatrix} \mid xw - yz = 1, x, y, z, w \in \mathbb{C} \right\}$. Then $\mathbb{C}[G] = \mathbb{C}[x, y, z, w] / (xw - yz - 1)$. 

    Let $f \in \mathbb{C}[G]$, then 
    $$
    \begin{aligned}
        X(f)(x, y, z, w) &= (x \partial_y + z \partial_w) f(x, y, z, w), \\
        Y(f)(x, y, z, w) &= (y \partial_x + w \partial_z) f(x, y, z, w), \\
        H(f)(x, y, z, w) &= (x \partial_x + z \partial_z - y \partial_y - w \partial_w) f(x, y, z, w).
    \end{aligned}
    $$

    Let $\rho$ be the adjoint representation of $\mathfrak{sl}_2(\mathbb{C})$. Under the basis $\{X, Y, H\}$, we have 
    $$
    \rho(X) = \begin{bmatrix}
        0 & 0 & -2 \\ 
        0 & 0 & 0 \\ 
        0 & 1 & 0 
    \end{bmatrix}, \quad 
    \rho(Y) = \begin{bmatrix}
        0 & 0 & 0 \\ 
        0 & 0 & 2 \\ 
        -1 & 0 & 0 
    \end{bmatrix}, \quad 
    \rho(H) = \begin{bmatrix}
        2 & 0 & 0 \\ 
        0 & -2 & 0 \\ 
        0 & 0 & 0 
    \end{bmatrix}
    $$ 
    and 
    $$
    \Tilde{\rho}\left(\begin{bmatrix}
        x & y \\ 
        z & w 
    \end{bmatrix}\right) = \begin{bmatrix}
        x^2 & -y^2 & -2xy \\ 
        -z^2 & w^2 & 2zw \\ 
        -xz & yw & xw + yz 
    \end{bmatrix}.
    $$

    Let $\mathbb{X}:\mathrm{SL}_2(\mathbb{C}) \rightarrow M_3(\mathbb{C})$ be a smooth solution to the following partial differential equations to obtain the gauge equivalence between $\rho$ and the trivial representation:
    $$
    \begin{aligned}
        &(x \partial_y + z \partial_w) \mathbb{X}(x, y, z, w) = \mathbb{X}(x, y, z, w) \rho(X), \\
        &(y \partial_x + w \partial_z) \mathbb{X}(x, y, z, w) = \mathbb{X}(x, y, z, w) \rho(Y), \\
        &(x \partial_x + z \partial_z - y \partial_y - w \partial_w) \mathbb{X}(x, y, z, w) = \mathbb{X}(x, y, z, w) \rho(H).
    \end{aligned}
    $$

    Then we have 
    $$
    \mathbb{X}(x, y, z, w) = \mathbb{X}(1, 0, 0, 1) \begin{bmatrix}
        x^2 & -y^2 & -2xy \\ 
        -z^2 & w^2 & 2zw \\ 
        -xz & yw & xw + yz 
    \end{bmatrix}
    $$ 
    and $\mathbb{X} \in M_n(\mathbb{C}[\mathrm{SL}_2(\mathbb{C})])$.
\end{example}

We know that every connected linear algebraic group is a semidirect product of a reductive group and a unipotent group. After the series of above examples of invariant D-modules on algebraic groups, we propose the following conjecture to conclude this section.

\begin{conjecture}
    Let $G$ be a connected linear algebraic group. We have
    $$
    \frac{\{\text{Invariant D-modules of rank } n \text{ over } G\}}{\sim} \cong \frac{\{\text{Invariant D-modules of rank } n \text{ over } \frac{G}{[G,G]}\}}{\sim}.
    $$

    Where $\sim$ denotes the relation of $D$-module isomorphisms. The above isomorphism is through the pull-back of representations via the quotient map $G \twoheadrightarrow G/[G,G]$.
\end{conjecture}

\section*{Acknowledgement}
I would like to thank Kifung Chan, Yujie Fu, Chenyu Bai and Zhongkai Tao for helpful discussions.

\bibliographystyle{plain}

\end{document}